\documentclass[a4paper,10pt]{amsart}

\usepackage{amssymb,amsthm,amsmath,amscd,amsfonts,bbm,mathrsfs}
\usepackage{stmaryrd}
\usepackage{hyperref}
\usepackage{xcolor}

\usepackage[cmtip,all]{xy}

\theoremstyle{plain}
\newtheorem{LeABC}{Lemma}
\newtheorem{ThmABC}[LeABC]{Theorem}
\newtheorem{Lemma}[equation]{Lemma}
\newtheorem{Thm}[equation]{Theorem}
\newtheorem{Prop}[equation]{Proposition}
\newtheorem{Cor}[equation]{Corollary}
\theoremstyle{definition}
\newtheorem{Defn}[equation]{Definition}
\newtheorem{Ass}[equation]{Assumption}
\newtheorem{Not}[equation]{Notation}
\newtheorem*{Not*}{Notation}
\theoremstyle{remark}
\newtheorem{Remark}[equation]{Remark}
\newtheorem*{Remark*}{Remark}
\newtheorem{Example}[equation]{Example}

\newcommand{\bLe}{\begin{Lemma}}
\newcommand{\eLe}{\end{Lemma}}
\newcommand{\bTh}{\begin{Thm}}
\newcommand{\eTh}{\end{Thm}}
\newcommand{\bPr}{\begin{Prop}}
\newcommand{\ePr}{\end{Prop}}
\newcommand{\bCo}{\begin{Cor}}
\newcommand{\eCo}{\end{Cor}}
\newcommand{\bDe}{\begin{Defn}}
\newcommand{\eDe}{\end{Defn}}
\newcommand{\bAs}{\begin{Ass}}
\newcommand{\eAs}{\end{Ass}}
\newcommand{\bNo}{\begin{Not}}
\newcommand{\eNo}{\end{Not}}
\newcommand{\bNoX}{\begin{Not*}}
\newcommand{\eNoX}{\end{Not*}}
\newcommand{\bReX}{\begin{Remark*}}
\newcommand{\eReX}{\end{Remark*}}
\newcommand{\bRe}{\begin{Remark}}
\newcommand{\eRe}{\end{Remark}}
\newcommand{\bEx}{\begin{Example}}
\newcommand{\eEx}{\end{Example}}
\newcommand{\bproof}{\begin{proof}}
\newcommand{\eproof}{\end{proof}}
\newcommand{\beqn}{\begin{equation}}
\newcommand{\eeqn}{\end{equation}}

\newcommand{\bToDo}{\begin{ToDo}\color{blue}}
\newcommand{\eToDo}{\end{ToDo}}

\newcommand{\bdone}{\bgroup\color{gray}}
\newcommand{\edone}{\egroup}
\newcommand{\btodo}{\bgroup\color{blue}}
\newcommand{\etodo}{\egroup}
\newcommand{\bnew}{\bgroup\color{violet}}
\newcommand{\enew}{\egroup}
\newcommand{\bnnew}{\bgroup\color{red}}
\newcommand{\ennew}{\egroup}
\newcommand{\bnnnew}{\bgroup\color{cyan}}
\newcommand{\ennnew}{\egroup}

\numberwithin{equation}{section}

\newcommand{\AAA}{\mathcal{A}}

\newcommand{\LLL}{\mathcal{L}}
\newcommand{\MMM}{\mathcal{M}}
\newcommand{\NNN}{\mathcal{N}}
\newcommand{\OOO}{\mathcal{O}}

\newcommand{\SSS}{\mathcal{S}}
\newcommand{\TTT}{\mathcal{T}}
\newcommand{\UUU}{\mathcal{U}}

\newcommand{\XXX}{\mathcal{X}}
\newcommand{\YYY}{\mathcal{Y}}
\newcommand{\ZZZ}{\mathcal{Z}}

\newcommand{\Fp}{\mathfrak{p}}
\newcommand{\Fq}{\mathfrak{q}}

\newcommand{\ZZ}{{\mathbb{Z}}}

\newcommand{\One}{{\mathbbm{1}}}

\newcommand{\pperf}{{\operatorname{-perf}}}
\newcommand{\pproj}{{\operatorname{-proj}}}

\DeclareMathOperator{\colim}{colim}
\DeclareMathOperator{\cone}{cone}
\DeclareMathOperator{\even}{even}

\DeclareMathOperator{\id}{id}

\DeclareMathOperator{\liset}{lis-\acute et}

\DeclareMathOperator{\modcat}{mod}

\DeclareMathOperator{\perf}{perf}
\DeclareMathOperator{\proj}{proj}
\DeclareMathOperator{\qc}{qc}
\DeclareMathOperator{\rep}{rep}

\DeclareMathOperator{\red}{red}
\DeclareMathOperator{\rel}{rel}

\DeclareMathOperator{\supp}{supp}
\DeclareMathOperator{\tr}{tr}

\DeclareMathOperator{\Ann}{Ann}
\DeclareMathOperator{\Aut}{Aut}

\DeclareMathOperator{\Comod}{Comod}
\DeclareMathOperator{\End}{End}

\DeclareMathOperator{\GL}{GL}
\DeclareMathOperator{\Hom}{Hom}

\DeclareMathOperator{\LF}{LF}

\DeclareMathOperator{\Mod}{Mod}

\DeclareMathOperator{\Perf}{Perf}
\DeclareMathOperator{\Qcoh}{Qcoh}
\DeclareMathOperator{\Rep}{Rep}

\DeclareMathOperator{\Sch}{Sch}
\DeclareMathOperator{\Spc}{Spc}

\DeclareMathOperator{\Spec}{Spec}

\begin{document}

\title
{Perfect complexes on finite flat affine groupoids}
\author{Eike Lau}
\date{\today}
\address{Eike Lau, Fakult\"{a}t f\"{u}r Mathematik,
Universit\"{a}t Bielefeld, D-33501 Bielefeld}

\begin{abstract}
We compute the Balmer spectrum of the category of perfect complexes on an algebraic stack admitting a finite locally free cover by an affine scheme and identify it with the homogeneous spectrum of the cohomology ring.
\end{abstract}

\maketitle

\setcounter{tocdepth}{1}
\tableofcontents

\section{Introduction}

For an algebraic stack $\XXX$ we consider the category $\Perf(\XXX)$ of perfect complexes on $\XXX$ as a tensor triangulated category.
Thick tensor ideals in $\Perf(\XXX)$ correspond to Thomason subsets of the Balmer spectrum $\Spc(\Perf(\XXX))$. 
The general theory 
gives a comparison map
\beqn
\label{Eq:rhoX-Intro}
\rho_\XXX\colon\Spc(\Perf(\XXX))\to\Spec^h(R_\XXX)
\eeqn
where $R_\XXX=H^*(\XXX,\OOO_\XXX)$ is the cohomology ring of $\XXX$ and $\Spec^h$ is the space of homogeneous prime ideals.
The aim of this note it to prove the following.

\begin{ThmABC}
\label{Th:Main}
If the algebraic stack $\XXX$ admits a finite locally free covering 
by an affine scheme, then $\rho_\XXX$ is a homeomorphism.
\end{ThmABC}

\subsubsection*{Algebraic translation}

The condition in Theorem \ref{Th:Main} means that $\XXX$ is presented by a groupoid of affine schemes with finite locally free source and target maps, or equivalently by a commutative Hopf algebroid $A\rightrightarrows B$ with finite locally free structure maps. In this situation, $\Perf(\XXX)$ is equivalent to $D^b(\LF(\XXX))$, the bounded derived category of the exact category of locally free $\OOO_\XXX$-modules of finite rank, and the latter correspond to comodules over $A\rightrightarrows B$ which are finite projective over $A$. This yields an algebraic formulation of Theorem \ref{Th:Main} without reference to algebraic stacks, but it will be important that equivalent groupoids have equivalent module categories.
As a special case, if $\XXX=[Y/G]$ is a quotient stack with $Y=\Spec A$ affine and $G$ a finite locally free group scheme over a common base ring, then $\LF(\XXX)$ is equivalent to the category of $G$-equivariant finite projective $A$-modules. 

\subsubsection*{Related results}

Theorem \ref{Th:Main} subsumes two recent results of a similar nature. If $\XXX$ is the quotient of an affine scheme by the action of a finite group, we recover the main result of \cite{Lau:Balmer}. If $\XXX$ is the classifying stack of a finite flat group scheme $G$ over a noetherian ring $k$, then $\Perf(\XXX)$ is equivalent to the category $\rep(G,k)$ of \cite{BBIKP:Lattices}, and the result is established in loc.\,cit.\ as a consequence of a stratification of a suitable ind-completion of $\rep(G,k)$. 

\subsubsection*{Concentrated stacks}

If $\XXX$ is a qcqs scheme, $\Spc(\Perf(\XXX))$ is homeomorphic to $|\XXX|$ by \cite{Thomason:Classification}. This extends to the case of a qcqs algebraic stack which is concentrated and satisfies the Thomason condition, by \cite[Theorem C]{Hall-Rydh:Telescope}.
The Thomason condition means that $D_{\qc}(\XXX)$ is compactly generated and every closed subset of $|\XXX|$ with quasi-compact complement is the support of a perfect complex; the stack $\XXX$ is concentrated if it has finite cohomological dimension for quasi-coherent sheaves, or equivalently if $\OOO_\XXX$ is compact in $D_{\qc}(\XXX)$. An algebraic stack $\XXX$ as in Theorem \ref{Th:Main} satisfies the Thomason condition by \cite[Theorem A]{Hall-Rydh:Perfect}, and $\XXX$ is concentrated iff $\XXX$ is tame, i.e.\ the (finite) stabiliser groups at geometric points are linearly reductive, by \cite[Theorem C]{Hall-Rydh:Algebraic}. The tame case of Theorem \ref{Th:Main} is also covered by \cite{Hall:Tame}.

\subsubsection*{Finite generation of cohomology.}

The following technical point will be important.
If the stack $\XXX$ in Theorem \ref{Th:Main} is of finite type over a noetherian ring $k$, which can be arranged by noetherian approximation, then $R_\XXX$ is a finitely generated $k$-algebra. Indeed, one can write $\XXX=[Z/{\GL_d}]$ with $Z=\Spec S$, and then $R_\XXX\cong H^*(\GL_d,S)$. 
If $\XXX$ is of finite type over $k$, the same holds for $S$, and the cohomology ring is finitely generated over $k$ iff it has bounded torsion by \cite{Kallen:Grosshans-family}. In our situation, the cohomology in positive degrees is annihilated by the degree of a finite locally free covering of $\XXX$ by an affine scheme.

\subsubsection*{Decomposition into fibers.}

The proof of Theorem \ref{Th:Main} follows an approach similar to \cite{Lau:Balmer}, by reduction to the known cases of an affine scheme \cite{Thomason:Classification} and of the classifying stack of a finite group scheme over a field \cite{BCR:Thick, Friedlander-Pevtsova:Pi-Supports}. 
The stack $\XXX$ has a coarse moduli space $X=\Spec C$ where $C=\Gamma(\XXX,\OOO_\XXX)$ is the ring of invariants of the Hopf algebroid $A\rightrightarrows B$; this is the degree zero part of the cohomolog ring $R_\XXX$. By functoriality we obtain a commutative diagram:
\beqn
\label{Eq:Intro-fund-diag}
\xymatrix@M+0.2em{
\Spc(\Perf(\XXX)) \ar[r]^-{\rho_\XXX} \ar[d] &
\Spec^h(R_\XXX) \ar[d] \\
\Spc(\Perf(X)) \ar[r]^-{\rho_X} &
\Spec(C)
}
\eeqn
Since the bottom arrow $\rho_X$ is a homeomorphism, we are led to consider the vertical fibers in \eqref{Eq:Intro-fund-diag}
over each $\Fp\in\Spec(C)$; these fibers will be denoted by a subscript $\Fp$.

\subsubsection*{Description of the fibers}

Let $\XXX(\Fp)$ be the reduced fiber of $\XXX$ over $\Fp$. This is the residual gerbe of $\XXX$ at the unique point of $\XXX$ over $\Fp$, in particular $\XXX(\Fp)$ is a gerbe with finite band over a field. The behavior of the vertical fibers in \eqref{Eq:Intro-fund-diag} is determined by the following facts which extend similar results of \cite{Lau:Balmer}.

\begin{LeABC}
\label{Le:Intro-fiber-Spec-h}
The natural morphism $\XXX(\Fp)\to\XXX$ induces a homeomorphism
\beqn
\Spec^h(R_{\XXX(\Fp)})\to\Spec^h(R_\XXX)_\Fp.
\eeqn
\end{LeABC}

The noetherian case of Lemma \ref{Le:Intro-fiber-Spec-h} is proved by a cohomological calculation that goes back to \cite{Benson-Habegger}; this case will be sufficient for the application. The general case follows by noetherian approximation. See \S\ref{Se:Fiber-cohomology}.

\begin{LeABC}
\label{Le:Intro-fiber-Spc}
If $\XXX$ is of finite type over a noetherian ring,
the natural morphism $\XXX(\Fp)\to\XXX$ induces a surjective map
\beqn
\label{Eq:Intro-fiber-Spc}
\Spc(\Perf(\XXX(\Fp)))\to\Spc(\Perf(\XXX))_\Fp.
\eeqn
\end{LeABC}

Lemma \ref{Le:Intro-fiber-Spc} is proved in \S\ref{Se:Fiber-Balmer} using a surjectivity criterion for a functor of tensor triangulated categories that detects tensor nilpotence of morphisms \cite{Balmer:Surjectivity}. 

\begin{LeABC}
\label{Le:Intro-rho-XXX-Fp}
The comparison map $\rho_{\XXX(\Fp)}$ is a homeomorphism.
\end{LeABC}

Lemma \ref{Le:Intro-rho-XXX-Fp} is proved in \S\ref{Se:Comparison-residual} by reduction to the case of a neutral gerbe, i.e.\ the classifying stack of a finite group scheme over a field, using a surjectivity criterion for a functor of tensor triangulated categories with a right adjoint \cite{Balmer:Surjectivity}.

Lemmas \ref{Le:Intro-fiber-Spec-h}, \ref{Le:Intro-fiber-Spc}, and \ref{Le:Intro-rho-XXX-Fp} show that $\rho_\XXX$ is bijective in the noetherian case, thus a homeomorphism. A posteriori, it follows that \eqref{Eq:Intro-fiber-Spc} is a homeomorphism in general.

\subsection*{Acknowledgements}
This work was 
funded by the Deutsche Forschungsgemeinschaft (DFG, German Research
Foundation) -- Project-ID 491392403 -- TRR 358.

\section{Preliminaries}

Let us fix some terminology and notation, along with some basic facts.

\subsection{Tensor triangular geometry, comparison map}

Let $\TTT$ be a tensor triangulated category and $R=\End^*(\One_\TTT)=\bigoplus_n\Hom(\One_\TTT,\One_\TTT[n])$ as a graded-commu\-ta\-tive ring. We consider the Balmer spectrum $\Spc(\TTT)$ defined in \cite{Balmer:Spectrum} and the comparison map $\rho_\TTT\colon\Spc(\TTT)\to\Spec^h(R)$ constructed in \cite{Balmer:Spectra3}, denoted $\rho^\bullet$ in loc.\,cit.
The category $\TTT$ is called noetherian if $\End^*(M)$ is a noetherian $R$-module for every $M\in\TTT$, or equivalently if the ring $R$ is noetherian and $\End^*(M)$ is finitely generated over $R$ for every $M$. If $\TTT$ is rigid and noetherian, the map $\rho_\TTT$ is a homeomorphism iff it is bijective by \cite[Corollary 2.8]{Lau:Balmer}.

\subsection{Sheaves of modules}

For a ringed site $\SSS$, let $\Qcoh(\SSS)$ denote the category of quasi-coherent $\OOO_{\SSS}$-modules \cite[\href{https://stacks.math.columbia.edu/tag/03DL}{Tag 03DL}]{stacks-project}, $D(\SSS)$ the derived category of $\OOO_\SSS$-modules, $D_{\qc}(\SSS)$ the full subcategory of $D(\SSS)$ consisting of complexes with quasi-coherent cohomology, and $\Perf(\SSS)$ the full subcategory of $D(\SSS)$ of perfect complexes \cite[\href{https://stacks.math.columbia.edu/tag/08G5}{Tag 08G5}]{stacks-project}. Then $\Perf(\SSS)\subseteq D_{\qc}(\SSS)$, moreover $\Perf(\XXX)$ is a rigid tensor triangulated category.
Perfect complexes are precisely the dualisable objects of $D(\SSS)$ by \cite[\href{https://stacks.math.columbia.edu/tag/0FPS}{Tag 0FPS}]{stacks-project}, or \cite[Lemma 4.3]{Hall-Rydh:Perfect} in the case of algebraic stacks.

\subsection{Algebraic stacks}

We assume that algebraic stacks are quasi-separated as in \cite{Laumon-MB:Champs}. This holds in our examples, and it simplifies the discussion of residual gerbes. For an algebraic stack $\XXX$ let $\XXX_{\liset}$ be the lisse-\'etale site \cite[\S12]{Laumon-MB:Champs} and $\Qcoh(\XXX)=\Qcoh(\XXX_{\liset})$ etc., in particular $\Perf(\XXX)=\Perf(\XXX_{\liset})$. The category $\Qcoh(\XXX)$ is Grothendieck abelian \cite[\href{https://stacks.math.columbia.edu/tag/0781}{Tag 0781}]{stacks-project} and therefore has enough injectives; cf.\ Remark \ref{Re:Qcoh-injectives}. In our examples, the lisse-\'etale site can be ignored by taking \eqref{Eq:PerfX-DQcohX-perf} below as a definition.

\bNoX
For an algebraic stack $\XXX$ we write $R_\XXX=\End^*(\OOO_\XXX)=H^*(\XXX,\OOO_\XXX)$ and
$\rho_\XXX=\rho_{\,\Perf(\XXX)}\colon\Spc(\Perf(\XXX))\to\Spec^h(R_\XXX)$ as in \eqref{Eq:rhoX-Intro}.
\eNoX

\section{Affine groupoids}

\label{Se:Affine-groupoids}

Let $\Sch$ be the category of schemes.
We note the following consequence of Artin's representability theorem for fppf groupoids \cite[Theorem 6.1]{Artin:Versal}. 

\bLe
\label{Le:aff-pres-stack}
For a fibered category $\XXX\to\Sch$ the following are equivalent.
\begin{enumerate}
\item
$\XXX$ is a quasi-compact algebraic stack with affine diagonal,
\item
$\XXX$ is the \'etale stack associated to a smooth groupoid of affine schemes,
\item
$\XXX$ is the fppf stack associated to an fppf groupoid of affine schemes.
\end{enumerate}
\eLe

When the equivalent conditions of Lemma \ref{Le:aff-pres-stack} hold, $\XXX$ will be called an affine-presented algebraic stack. These stacks are called geometric in \cite{Lurie:Tannaka}. In this section we recall the description of quasi-coherent sheaves and their cohomology and of perfect complexes on such stacks.

\subsection{Quasi-coherent sheaves and comodules}
\label{Se:Qcoh}

Let $\XXX$ be an affine-pre\-sented algebraic stack, let $\pi\colon Y_0\to\XXX$ be an fppf covering where $Y_0$ is an affine scheme, and let $Y_1=Y_0\times_{\XXX}Y_0$. Then $Y_1\rightrightarrows Y_0$ is an fppf groupoid of affine schemes whose fppf quotient $[Y_0/Y_1]$ is isomorphic to $\XXX$. The groupoid $Y_1\rightrightarrows Y_0$ corresponds to a commutative fppf Hopf algebroid $A\rightrightarrows B$ via $Y_0=\Spec A$ and $Y_1=\Spec B$.

The category $\Qcoh(\XXX)$ is equivalent to the category $\Qcoh(Y_1\rightrightarrows Y_0)$ of quasi-coherent sheaves on the groupoid as in \cite[\href{https://stacks.math.columbia.edu/tag/0441}{Tag 0441}]{stacks-project}, or equivalently the category $\Comod(A\rightrightarrows B)$ of comodules over the Hopf algebroid.
Under this identification, the pull-back functor $\pi^*\colon\Qcoh(\XXX)\to\Qcoh(Y_0)$ corresponds to the forgetful functor $\Comod(A\rightrightarrows B)\to\Mod(A)$. It has a right adjoint $\pi_*\colon\Qcoh(Y_0)\to\Qcoh(\XXX)$, corresponding to the functor $\Mod(A)\to\Comod(A\rightrightarrows B)$, $M\mapsto M\otimes_AB$.

Let $\LF(\XXX)$ denote the category of locally free $\OOO_{\XXX}$-modules.
Then $\LF(\XXX)\subseteq\Qcoh(\XXX)$, and $F\in\Qcoh(\XXX)$ is locally free iff $\pi^*F$ is locally free. Thus $\LF(\XXX)$ is equivalent to the category $\Comod(A\rightrightarrows B)_{A\pproj}$ of comodules which are finite projective over $A$.

\bRe
\label{Re:Qcoh-injectives}
One sees directly that the category $\Qcoh(\XXX)$ has enough injectives, and these are direct summands of $\pi_*I$ with injective $I\in\Qcoh(Y_0)$. Indeed,
the functor $\pi_*$ of quasi-coherent modules preserves injectives since $\pi^*$ is exact. For $F\in\Qcoh(\XXX)$, the natural map $F\to\pi_*\pi^*F$ is a monomorphism since $\pi$ is affine and faithfully flat. Choose a monomorphism $\pi^*F\to I$ with injective $I$; then $F\to\pi_*I$ is a monomorphism with $\pi_*I$ injective.
\eRe

\subsection{Derived categories of sheaves}
\label{Se:derived-cat-sheaves}

Let $\XXX$ be an affine-presented algebraic stack and $\pi\colon Y_0\to\XXX$ an fppf covering by an affine scheme as above.
By \cite[Theorem 3.8]{Lurie:Tannaka} or \cite[Theorem C.1]{Hall-Neeman-Rydh} there is an equivalence
\beqn
\label{Eq:D+Qcoh-D+qc}
D^+(\Qcoh(\XXX))\cong D^+_{\qc}(\XXX),
\eeqn
which induces an equivalence $D^b(\Qcoh(\XXX))\cong D^b_{\qc}(\XXX)$. Let 
\beqn
\label{Eq:DQcohXperf}
D(\Qcoh(\XXX))_{\perf}\subseteq D(\Qcoh(\XXX))
\eeqn
be the full subcategory of all complexes $C$ such that $\pi^*C$ is perfect.
This is independent of the choice of $\pi$ since perfectness over schemes is fppf local by \cite[\href{https://stacks.math.columbia.edu/tag/068T}{Tag 068T}]{stacks-project}. 
We have $D(\Qcoh(\XXX))_{\perf}\subseteq D^b(\Qcoh(\XXX))$, and \eqref{Eq:D+Qcoh-D+qc} induces an equivalence
\beqn
\label{Eq:PerfX-DQcohX-perf}
D(\Qcoh(\XXX))_{\perf}\cong\Perf(\XXX).
\eeqn
This yields a description of $\Perf(\XXX)$ independent of the lisse-\'etale site.

\subsection{Cohomology of quasi-coherent sheaves}

Let $\XXX$ be again an affine-pre\-sent\-ed algebraic stack.
For $F\in\Qcoh(\XXX)$, the sheaf cohomology groups $H^n(\XXX,F)$ form a derived functor on the category $\Qcoh(\XXX)$ by \eqref{Eq:D+Qcoh-D+qc}, and they can be computed by \v Cech cohomology as follows. 
For an fppf covering $\pi\colon Y_0\to\XXX$ with affine $Y_0$
let $Y_\bullet$ be the \v Cech nerve of $\pi$ and
\[
\check H^n(Y_0\to\XXX,F)=H^n(\Gamma(Y_\bullet,F)).
\]

\bLe
\label{Le:Cech}
For $F\in\Qcoh(\XXX)$ there is a natural isomorphism 
\[
H^n(\XXX,F)\cong\check H^n(Y_0\to\XXX,F).
\]
\eLe

\bproof
This follows from \cite[Corollary 9.2.4]{Olsson:Stacks}. One can also argue directly as follows.
If $f\colon Z\to\XXX$ is flat with affine $Z$, for $N\in\Qcoh(Z)$ we have $H^n(\XXX,f_*N)=0$ for $n\ge 1$. Indeed, the functor $f_*\colon\Qcoh(Z)\to\Qcoh(\XXX)$ is exact and preserves injective objects since $f^*$ is exact, and the functor $\Gamma(\XXX,f_*-)=\Gamma(Z,-)$ is exact.
Let $\pi_\bullet\colon Y_\bullet\to\XXX$ be the natural morphism. Since each $Y_n$ is affine and $\pi_n$ is flat, $0\to F\to\pi_{\bullet*}\pi_\bullet^*F$ is an acyclic resolution, 
and $\Gamma(\XXX,\pi_{\bullet*}\pi_\bullet^*F)=\Gamma(Y_\bullet,F)$.
\eproof

\bCo
\label{Co:H^nXf_*M}
For any morphism $f\colon Z\to\XXX$ with affine $Z$ and for $N\in\Qcoh(Z)$ we have $H^n(\XXX,f_*N)=0$ for $n\ge 1$. 
\eCo

\bproof
The group $\check H^n(Y_0\to\XXX,f_*N)$ identifies with $\check H^n(W_0\to Z,N)$ where $W_0=Y_0\times_\XXX Z$ is affine and flat over $Z$.
\eproof

\subsection{Quotient case}
\label{Se:[Y/G]-gen}

Let $k$ be a ring, $Y=\Spec A$ an affine $k$-scheme, and $G$ a flat affine group scheme of finite presentation over $k$ that acts on $Y$. The quotient $\XXX=[Y/G]$ is an affine-presented algebraic stack with fppf covering $\pi\colon Y\to\XXX$. 
In this situation, 
the associated categories are as follows.

Let $\Rep_A(G)$ be the category of $G$-equivariant $A$-modules and $\Rep^{\proj}_A(G)$ the category of $G$-equivariant finite projective $A$-modules. Then 
$\Qcoh(\XXX)$ is equivalent to $\Rep_A(G)$ via $F\mapsto M=\Gamma(Y,\pi^*F)$,
this restricts to an equivalence between $\LF(\XXX)$ and $\Rep_A^{\proj}(G)$, and $H^n(\XXX,F)$ is isomorphic to the algebraic group cohomology $H^n(G,M)$ as defined in \cite{Jantzen:Book}.
Moreover, \eqref{Eq:PerfX-DQcohX-perf} reads
\beqn
\label{Eq:PerfX-DOepAG-Aperf}
\Perf(\XXX)\cong D(\Rep_A(G))_{A\pperf}
\eeqn
with $D(\Rep_A(G))_{A\pperf}$ denoting the full subcategory of the derived category of $\Rep_A(G)$ where the underlying complex of $A$-modules is perfect.

As a special case, for $Y=\Spec k$, with trivial $G$-action, $[Y/G]=BG$ is the classifying stack of $G$, and $\Rep_k(G)$ is the category of representations of $G$ over $k$.

\subsection{Resolution property}
\label{Se:resolution}

An algebraic stack $\XXX$ has the resolution property if every quasi-coherent sheaf of finite type is a quotient of a locally free sheaf.
Assume that $\XXX$ is qcqs with affine stabiliser groups at closed points. 
By \cite[Theorem A]{Gross:Tensor-generators}, such a stack $\XXX$ has the resolution property iff $\XXX$ is basic, meaning that $\XXX\cong[U/{\GL_d}]$ where $U$ is a quasi-affine scheme. This property implies that $\XXX$ has affine diagonal, thus $\XXX$ is affine-presented in the above sense.

\bLe
\label{Le:resolution-Perf-DbLF}
If $\XXX$ is an affine-presented algebraic stack with the resolution property, then the natural functor $D^b(\LF(\XXX))\to\Perf(\XXX)$ is an equivalence;
here $D^b$ is the bounded derived category of an exact category as in \cite[\S11]{Keller:Derived} or \cite[\S4.1]{Krause:Homological}.
\eLe

\bproof
This is proved in \cite[\href{https://stacks.math.columbia.edu/tag/0F8I}{Tag 0F8I}]{stacks-project} when $\XXX$ is a scheme, but the proof carries over to the stack case, using that $\Perf(\XXX)\cong D(\Qcoh(\XXX))_{\perf}$ by \eqref{Eq:PerfX-DQcohX-perf} and that every quasi-coherent sheaf on $\XXX$ is the filtered union of its quasi-coherent subsheaves of finite type by \cite{Rydh:Approximation-sheaves}.
\eproof

\section{Finite locally free affine groupoids}

In this section we recall basic properties of algebraic stacks presented by finite locally free groupoids of affine schemes.

\subsection{Well-nigh affine stacks}
\label{Se:wna-stacks}

Following \cite[\href{https://stacks.math.columbia.edu/tag/0DUL}{Tag 0DUL}]{stacks-project}, an algebraic stack $\XXX$ will be called well-nigh affine if there is a surjective finite locally free morphism $Y_0\to\XXX$ where $Y_0$ is an affine scheme.\footnote{As noted in loc.\,cit.\ this terminology reflects that this class of stacks is rarely central.}
This holds iff $\XXX$ is presented by a finite locally free groupoid of affine schemes $Y_1\rightrightarrows Y_0$, given by $Y_1=Y_0\times_\XXX Y_0$; in particular $\XXX$ is affine-presented as in \S \ref{Se:Affine-groupoids}. Such groupoids are equivalent to commutative finite locally free Hopf algebroids $A\rightrightarrows B$.

\bRe
\label{Re:degree}
One can assume that $Y_0\to\XXX$ as above has constant degree because in general, the degree is locally constant on $\XXX$, and it can be adjusted by passing to suitable multiples of $Y_0$ over each locus where the degree is constant.
\eRe

\bEx
\label{Ex:[Y/G]-wna}
Let $k$ be a ring, $Y$ an affine $k$-scheme, and $G$ a finite locally free group scheme over $k$ that acts on $X$. Then $[Y/G]$ is a well-nigh affine stack with finite locally free covering $Y\to[Y/G]$.
\eEx

\subsection{Coarse moduli space}
\label{Se:Coarse}

Let $\XXX$ be a well-nigh affine stack, presented by a finite locally free groupoid of affine schemes $Y_1\rightrightarrows Y_0$, corresponding to a commutative finite locally free Hopf algebroid $A\rightrightarrows B$ as in \S\ref{Se:wna-stacks}.
We consider the ring of global sections
\beqn
\label{Eq:C-Gamma-ker}
C=\Gamma(\XXX,\OOO_{\XXX})=\ker(A\rightrightarrows B)
\eeqn
and the affine scheme $X=\Spec C$. The natural morphism $f\colon\XXX\to X$ is a uniform categorical moduli space in the category of algebraic spaces; moreover $f$ is separated, quasi-compact, and a universal homeomorphism; see \cite[\href{https://stacks.math.columbia.edu/tag/0DUR}{Tag 0DUR}]{stacks-project}. 
In more detail, the homomorphism $C\to A$ is integral, the diagram of topological spaces $|Y_1|\rightrightarrows|Y_0|\to|X|$ is a coequalizer, and for any ring homomorphism $C\to C'$ and $C^1=\ker(A\otimes_CC'\rightrightarrows B\otimes_CC')$, the ring homomorphism $C'\to C^1$ induces a universal homeomorphism on spectra; see the proof of \cite[\href{https://stacks.math.columbia.edu/tag/0DUP}{Tag 0DUP}]{stacks-project}. If $C\to C'$ is flat, then $C'=C^1$. 
If the groupoid $Y_1\rightrightarrows Y_0$ has trivial stabilisers, then $\XXX$ is an algebraic space and hence $\XXX\cong X$; see \cite[\href{https://stacks.math.columbia.edu/tag/03BM}{Tag 03BM}]{stacks-project} for a direct proof. In other words, a well-nigh affine algebraic space is an affine scheme.

The original reference for these facts is \cite[Expos\'e V]{SGA3}.

\bRe
\label{Re:XXX-k-finite-type}
The stack $\XXX$ is of finite type over a noetherian ring $k$ iff $A\rightrightarrows B$ is a Hopf algebroid in the category of $k$-algebras of finite type. In this case, $C$ is
of finite type over $k$ as well, since the composition
$k\to C\to A$ is of finite type and $C\to A$ is integral, which implies that $k\to C$ is of finite type.
\eRe

\subsection{Quotient presentation}

\bLe
\label{Le:XXX=Z/GLd}
Every well-nigh affine stack $\XXX$ can be written as $\XXX\cong[Z/{\GL_d}]$ where $Z$ is an affine scheme.
\eLe

\bproof
Assume that $\pi\colon Y\to\XXX$ is surjective finite locally free of constant rank $d$ where $Y$ is affine; see Remark \ref{Re:degree}. Let $Z\to\XXX$ be the $\GL_d$-torsor that parametrizes trivializations of the locally free $\OOO_{\XXX}$-module $\pi_*\OOO_{Y}$. The inertia of $\XXX$ acts freely, thus faithfully on the fibers of $Y\to\XXX$, so it acts freely on the fibers of $Z\to\XXX$. Hence $Z$ has trivial inertia.
Let $W=Y\times_\XXX Z$. Then $W\to Y$ is affine, hence $W$ is affine. Now $W\to Z$ is surjective finite locally free, hence $Z$ is a well-nigh affine algebraic space, thus an affine scheme.
\eproof

\subsection{Noetherian approximation}

\bDe
\label{De:limit-stack-affine}
Let $\YYY$ be an algebraic stack, $g\colon\ZZZ\to\YYY$ an affine morphism, and $g_i\colon\ZZZ_i\to\YYY$ a cofiltered system of affine morphisms. We write $\ZZZ\cong\lim_i\ZZZ_i$ if $g_*\OOO_\ZZZ\cong\colim_i g_{i*}\OOO_{\ZZZ_i}$ as quasi-coherent $\OOO_{\YYY}$-algebras.
\eDe

\bLe
\label{Le:approximation}
Let $\XXX$ be a well-nigh affine stack and $\pi\colon Y\to\XXX$ surjective finite locally free with affine $Y$. There are a well-nigh affine stack $\XXX_0$ of finite presentation over $\ZZ$, a surjective finite locally free morphism $\pi\colon Y_0\to\XXX_0$ with affine $Y_0$, and an affine morphism $\XXX\to\XXX_0$ such that $Y\cong\XXX\times_{\XXX_0}Y_0$. Moreover, $\XXX\to\XXX_0$ is a cofiltered limit of affine morphisms of finite presentation $\XXX_i\to\XXX_0$ as in Definition \ref{De:limit-stack-affine}.
\eLe

\bproof
This is a consequence of the general approximation results of \cite{Rydh:Noeth-approx}. In the present situation, one can argue directly as follows.

Let $\XXX=[Z/{\GL_d}]$ with $Z=\Spec S$; see Lemma \ref{Le:XXX=Z/GLd}. The ring $S$ is the filtered colimit of its $\GL_d$-invariant subrings of finite type $S_i$, since $S$ is the filtered colimit of the invariant finite submodules by \cite[Proposition 15.4]{Laumon-MB:Champs}.
Let $Z_i=\Spec S_i$ and $\XXX_i=[Z_i/{\GL_d}]$. Then $\XXX=\lim_i\XXX_i$.

The surjective finite locally free morphism $Y\to\XXX$ corresponds to a $\GL_d$-equivariant surjective finite locally free morphism $W\to Z$ via $W=Y\times_\XXX Z$ and $Y=[W/{\GL_d}]$. The morphism $W\to Z$ comes from a $\GL_d$-equivariant surjective finite locally free morphism $W_i\to Z_i$ for some $i$, which gives $W_j\to Z_j$ for $j\ge i$. Let $Y_j=[W_j/{\GL_d}]$ as an algebraic stack. Then $Y_j\to\XXX_j$ is surjective finite locally free.
We have to show that $Y_j$ is affine for some $j$.

Let $Y=\Spec A$ and $A=\colim_m A_m$ as a filtered colimit of rings of finite type, $Y'_m=\Spec A_m$. The $\GL_d$-torsor $W\to Y$ comes from a $\GL_d$-torsor $W'_{m}\to Y'_{m}$ for some $m$. 
The morphism $W\to W'_m$ factors over $W_j\to W'_m$ for some $j$, which gives an affine morphism $Y_j\to Y'_m$. Thus $Y_j$ is affine.
\eproof

\bRe
\label{Re:approximation-LF}
For a filtered limit $\XXX=\lim_i\XXX_i$ as in Lemma \ref{Le:approximation}, the natural functor $\colim_i\LF(\XXX_i)\to\LF(\XXX)$ is an equivalence by the comodule description in \S\ref{Se:Qcoh}.
\eRe

\section{Projective and locally free modules}

Every well-nigh affine stack $\XXX$ has the resolution property as a consequence of Lemma \ref{Le:XXX=Z/GLd} (see \S\ref{Se:resolution}), hence $\Perf(\XXX)$ is equivalent to $D^b(\LF(\XXX))$ by Lemma~\ref{Le:resolution-Perf-DbLF}. 
In this section we prove these properties directly, starting from the observation that there are enough projective quasi-coherent modules, which also yields a computation of homomorphisms in $\Perf(\XXX)$ (Remark \ref{Re:PQ-special}).

\subsection{Finite duality}

Let us recall the duality for a finite locally free morphism.

\bLe
\label{Le:finite-duality}
Let $f\colon\YYY\to\XXX$ be a finite locally free morphism of algebraic stacks and let $f^*\colon\Qcoh(\XXX)\to\Qcoh(\YYY)$ be the pullback functor.
\begin{enumerate}
\item
There is a sequence of adjoint functors $f^*\dashv f_*\dashv f^!$.
\item
The functors $f^*$, $f_*$, $f^!$ are exact and preserve modules of finite type.
\item
There is a natural isomorphism $f^!\MMM\cong f^*\MMM\otimes_{\OOO_X}f^!\OOO_{\XXX}$.
\item
If $f'\colon\YYY'\to\XXX'$ is the base change of $f$ under a morphism of algebraic stacks $g\colon\XXX'\to\XXX$, and $h\colon\YYY'\to\YYY$ is the projection, then there are natural base change isomorphisms $g^*f_*\cong f'_*h^*$ and $h^*f^!\cong f'^!g^*$.
\end{enumerate}
\eLe

\bproof
If $\XXX=\Spec A$ is affine, thus $\YYY=\Spec B$, under the identifications $\Qcoh(\XXX)=\Mod(A)$ and $\Qcoh(\YYY)=\Mod(B)$ we have $f^*(M)=M\otimes_AB$, $f_*(N)=N$ as an $A$-module, and $f^!(M)=\Hom_A(B,M)$ as a $B$-module. These functors in the affine case have the stated properties, including the base change isomorphisms when $\XXX'$ and thus $\YYY'$ are affine as well. The lemma follows by fppf descent.
\eproof

\bLe
\label{Le:counit!}
Let $f\colon\YYY\to\XXX$ be a surjective finite locally free morphism of algebraic stacks. For $F\in\Qcoh(\XXX)$ the counit $f_*f^!F\to F$ is surjective.
\eLe

\bproof
We can assume that $\XXX=\Spec A$ and $\YYY=\Spec B$. For an $A$-module $M$ we have to show that $\Hom_A(B,M)\to M$, $\varphi\mapsto\varphi(1)$ is surjective.
Since $f$ is surjective, $A\to B$ is faithfully flat, so the quotient module $B/A$ is $A$-flat (the sequence of $A$-modules $0\to A\to B\to B/A\to 0$ splits over $B$), thus $B/A$ is $A$-projective, and the assertion follows.
\eproof

\subsection{Projective modules}

Let $\XXX$ be an algebraic stack and $\pi\colon Y\to\XXX$ surjective finite locally free where $Y$ is affine, so $\XXX$ is well-nigh affine (\S\ref{Se:wna-stacks}).

\bLe 
\label{Le:projective-Qcoh-XXX}
We consider the functors $\pi_*$ and $\pi^*$ of quasi-coherent sheaves.
\begin{enumerate}
\item
$\pi_*$ and $\pi^*$ preserve projectives and projectives of finite type.
\item
The category $\Qcoh(\XXX)$ has enough projectives, and these are the direct summands of $\pi_*P$ for projective $P$. 
\item
Every $F\in\Qcoh(\XXX)$ of finite type is a quotient of a projective of finite type, and these are the direct summands of $\pi_*P$ where $P\in\Qcoh(Y)$ is projective of finite type.
\item
Projective objects of finite type in $\Qcoh(\XXX)$ are locally free.
\end{enumerate}
\eLe

\bproof
The functors $\pi_*$ and $\pi^*$ preserve modules of finite type, and they preserve projectives since the right adjoint functors $\pi^!$ and $\pi_*$ are exact; see Lemma \ref{Le:finite-duality}. For $F\in\Qcoh(\XXX)$ we choose a projective $P\in\Qcoh(Y)$ and a surjective homomorphism $P\to\pi^!F$. Then $\pi_*P\to\pi_*\pi^!\MMM\to\MMM$ is surjective by Lemma \ref{Le:counit!}. If $F$ is of finite type, the same holds for $\pi^!F$, so we can choose $P$ of finite type; then $Q=\pi_*P$ is of finite type. Moreover, this $Q$ is locally free, and the last assertion follows.
\eproof

\bLe
\label{Le:D-C}
For every complex $C$ in $\Qcoh(\XXX)$ such that $\pi^*C$ is perfect, there is a quasi-isomorphism $D\to C$ where $D$ is a bounded above complex of projective objects in $\Qcoh(\XXX)$ of finite type.
\eLe

\bproof
This is straightforward from Lemma \ref{Le:projective-Qcoh-XXX}.
Use induction over the width of the Tor-amplitude of $\pi^*C$. 
For width $0$ we have to show that every $M\in\LF(\XXX)$ has a resolution by projectives of finite type, which is clear by Lemma \ref{Le:projective-Qcoh-XXX}.
Assume that $\pi^*(C)$ has Tor-amplitude in $[-n,0]$ with $n>0$. Then $H^0(C)$ is of finite type. By Lemma \ref{Le:projective-Qcoh-XXX} 
there is  $Q\in\Qcoh(\XXX)$ projective of finite type and a homomorphism $Q\to C$ such that $Q\to H^0(C)$ is surjective. Let $C'=\cone(Q\to C)$. Then $\pi^*C'$ is perfect with Tor-amplitude in $[-n,-1]$. Since $C$ is homotopy equivalent to $\cone(C'[-1]\to Q)$, the assertion follows by induction.
\eproof

\subsection{Perfect complexes}

Let $\XXX$ be again a well-nigh affine stack.

\bPr
\label{Pr:DbLF-DQcoh-perf}
The natural functor $j\colon D^b(\LF(\XXX))\to\Perf(\XXX)$ is an equivalence.
\ePr

\bproof
Let $\Qcoh(\XXX)_{\proj}\subseteq\Qcoh(\XXX)$ be the full subcategory of projectives of finite type. The natural functor $K^-(\Qcoh(\XXX)_{\proj})\to D(\Qcoh(\XXX))$ is fully faithful. Let $K^-(\Qcoh(\XXX)_{\proj})_{\perf}$ be the inverse image of $D(\Qcoh(\XXX))_{\perf}$ under this functor.
By Lemma \ref{Le:D-C} and \eqref{Eq:PerfX-DQcohX-perf}
we obtain an equivalence
\[
K^-(\Qcoh(\XXX)_{\proj})_{\perf}\cong D(\Qcoh(\XXX))_{\perf}\cong\Perf(\XXX).
\]
Let us construct an inverse of $j$. Let $C\in K^-(\Qcoh(\XXX)_{\proj})_{\perf}$. For sufficiently small $n$, the truncation $\tau_{\ge n}C$ is perfect with locally free components in degrees $>n$, thus locally free components in all degrees, moreover $\tau_{\ge n}C$ stabilises in $D^b(\LF(\XXX))$ for decreasing $n$. The inverse of $j$ is given by $C\mapsto\lim_n\tau_{\ge n}C$.
\eproof

\bRe
\label{Re:PQ-special}
Homomorphisms in $\Perf(\XXX)$ can be computed in the bounded homotopy category $K^b(\Qcoh(\XXX))$ as follows. Let $P,Q\in\Perf(\XXX)$. Up to isomorphism we can assume that $P$ is a bounded above complex of projectives in $\Qcoh(\XXX)$ of finite type and that $Q$ is a bounded complex in $\LF(\XXX)$, by \eqref{Eq:PerfX-DQcohX-perf}, Lemma \ref{Le:D-C} and Proposition \ref{Pr:DbLF-DQcoh-perf}. 
If $Q$ lies in degrees $>n$, then
\beqn
\label{Eq:Hom-Perf-PQ}
\Hom_{\Perf(\XXX)}(P,Q)=\Hom_{K^b(\Qcoh(\XXX))}(\sigma_{\ge n}P,Q).
\eeqn
\eRe

\bPr
\label{Pr:PerfXXX-colim}
For a filtered limit $\XXX=\lim_i\XXX_i$ as in Lemma \ref{Le:approximation}, the natural functor $\colim_i\Perf(\XXX_i)\to\Perf(\XXX)$ is an equivalence.
\ePr

\bproof
Using Remark \ref{Re:approximation-LF}, the functor is essentially surjective by Proposition \ref{Pr:DbLF-DQcoh-perf} and fully faithful by Remark \ref{Re:PQ-special}. 
\eproof

\subsection{Locally free modules}

For completeness we record the following fact.

\bPr
For a well-nigh affine stack $\XXX$,
the exact category $\LF(\XXX)$ is a Frobenius category, i.e.\ it has enough injectives and projectives, and these coincide.
\ePr

\bproof
Let $\pi\colon Y\to\XXX$ be surjective finite locally free with affine $Y$. 
The functor $\pi_*$ restricts to a functor $\pi_*^{\LF}\colon\LF(Y)\to\LF(\XXX)$ with exact left adjoint given by the restriction of $\pi^*$. Hence $\pi_*^{\LF}$ preserves injectives, moreover $\pi_*^{\LF}$ preserves projectives since this holds for $\pi_*$. The exact category $\LF(Y)$ is split. Hence direct summands of $\pi_*P$ are projective and injective in $\LF(\XXX)$, and all projectives are of this kind by Lemma \ref{Le:projective-Qcoh-XXX}. Thus projective implies injective, and the duality involution gives the converse.
\eproof

\subsection{Gorenstein coverings}

The following remarks will not be used later, but they may provide useful illustrations.

\bRe
\label{Re:Comment-f!*}
Let $f\colon\YYY\to\XXX$ be a finite locally free morphism of algebraic stacks and $f^*\dashv f_*\dashv f^!$ as in Lemma \ref{Le:finite-duality}. The following are equivalent.
\begin{enumerate}
\item
The functor $f^!$ preserves locally free sheaves.
\item
The $\OOO_\YYY$-module $f^!\OOO_\XXX$ is invertible. 
\item 
The morphism $f$ is Gorenstein.
\end{enumerate}
Indeed, by the local flatness criterion this can be verified when $\XXX=\Spec k$ for a field $k$ and $\YYY=\Spec B$ for a finite $k$-algebra $B$, and the assertion is the classical fact that $B$ is Gorenstein iff $\Hom_k(B,k)$ is a free $B$-module, necessarily of rank $1$.

If $\LLL=f^!\OOO_\XXX$ is invertible, the sequence of adjoint functors extends infinitely on both sides: $\ldots\dashv f_!\dashv f^*\dashv f_*\dashv f^!\dashv \ldots$ with $f_!(\NNN)=f_*(\NNN\otimes\LLL)$ etc., and these functors preserve locally free modules.
\eRe

\bEx
Let $\XXX=[Y/G]$ as in Example \ref{Ex:[Y/G]-wna}, so $Y=\Spec A$. Then the morphism $\pi\colon Y\to\XXX$ is Gorenstein since flat group schemes of finite presentation are syntomic by \cite[Proposition 27.26]{GW:AG2}, so the functor $\pi_!$ of Remark \ref{Re:Comment-f!*} exists. Under the identification $\LF(\XXX)=\Rep_A^{\proj}(G)$ of \S\ref{Se:[Y/G]-gen}, the functor $\pi^*$ is the restriction $\Rep_A^{\proj}(G)\to\Mod^{\proj}(A)$, thus $\pi_*$ corresponds to induction and $\pi_!$ corresponds to coinduction as defined in \cite{Jantzen:Book}.
\eEx

\section{Finite generation of cohomology}

Let $\XXX$ be a well-nigh affine stack, and let $\pi\colon Y\to\XXX$ be surjective finite locally free of degree $d$ with affine $Y$; see Remark \ref{Re:degree}.

\bLe
\label{Le:HnXOX-d}
The group $H^n(\XXX,\OOO_\XXX)$ is annihilated by $d$ for $n\ge 1$.
\eLe

\bproof
Since $\pi_*\OOO_Y$ is a locally free $\OOO_\XXX$-algebra of rank $d$, there is a trace homomorphism $\tr\colon\pi_*\OOO_Y\to\OOO_\XXX$ where $\tr(c)$ is the trace of the endomorphism of $\pi_*\OOO_Y$ given by multiplication with $c$. The composition 
\beqn
\label{Eq:trace}
\OOO_\XXX\xrightarrow{\pi^*}\pi_*\OOO_Y\xrightarrow{\tr}\OOO_\XXX
\eeqn
is multiplication by $d$. Since $\pi_*\OOO_Y$ has trivial higher cohomology by Corollary \ref{Co:H^nXf_*M} (or by the proof of Lemma \ref{Le:Cech}), the lemma follows.
\eproof

\bEx
If $\XXX=BG$ for a finite locally free group scheme $G=\Spec\OOO_G$ over a ring $k$, then \eqref{Eq:trace} corresponds to the sequence of $G$-modules $k\to\OOO_G\xrightarrow\tau k$ where $\tau$ is the trace. 
\eEx

\bCo
\label{Co:HnXN-d}
For $F\in\Qcoh(\XXX)$, $H^n(\XXX,F)$ is annihilated by $d$ for $n\ge 1$.
\eCo

\bproof
Let $\AAA=\OOO_\XXX\oplus F$ as an $\OOO_\XXX$-algebra with $F^2=0$ and $\ZZZ=\Spec_\XXX(\AAA)$. Then $g\colon\ZZZ\to\XXX$ is affine, and $Y\times_{\XXX}\ZZZ\to\ZZZ$ is a finite locally free covering of degree $d$ by an affine scheme, in particular $\ZZZ$ is well-nigh affine. Hence $d$ annihilates  
$H^n(\ZZZ,\OOO_\ZZZ)=H^n(\XXX,g_*\OOO_\ZZZ)=H^n(\XXX,\OOO_\XXX)\oplus H^n(\XXX,F)$.
\eproof

\bPr
\label{Pr:R-finite}
If $\XXX$ is of finite type over a noetherian ring $k$, the $k$-algebra $R=H^*(\XXX,\OOO_\XXX)$ is finitely generated, for a coherent $\OOO_{\XXX}$-module $F$ the $R$-module $H^*(\XXX,F)$ is finite, 
and the category $\Perf(\XXX)$ is noetherian.
\ePr

\bproof
The $k$-algebra $R_0=H^0(\XXX,\OOO_\XXX)$ is of finite type by Remark \ref{Re:XXX-k-finite-type}, and $R$ is of bounded torsion by Lemma \ref{Le:HnXOX-d}. We can write $\XXX\cong[Z/{\GL_d}]$ with $Z=\Spec S$ by Lemma \ref{Le:XXX=Z/GLd}. Here $S$ is of finite type over $k$ since this holds for $\XXX$, and $R=H^*(\GL_d,S)$ by \S\ref{Se:[Y/G]-gen}. For this graded $k$-algebra, bounded torsion implies finite generation by \cite[Theorem 10.5]{Kallen:Grosshans-family}. In particular, $R$ is a noetherian graded ring.

For given $F\in \Qcoh(\XXX)$ of finite type 
let $\ZZZ=\Spec_\XXX(\OOO_\XXX\oplus F)$ as in the proof of Corollary \ref{Co:HnXN-d}. Then $\ZZZ\to\XXX$ is affine and of finite type, so $\ZZZ$ is well-nigh affine and of finite type over $k$, thus $H^*(\ZZZ,\OOO_\ZZZ)=R\oplus H^*(\XXX,F)$ is a noetherian ring, hence $H^*(\XXX,F)$ is finitely generated as an $R$-module.

For $M\in\Perf(\XXX)$, the $R$-module $\End^*(M)$ is finitely generated since by \eqref{Eq:PerfX-DQcohX-perf} and Proposition \ref{Pr:DbLF-DQcoh-perf}, $\End^*(M)$ is a successive extension of submodules and quotients of shifts of $H^*(\XXX,E^\vee\otimes F)$ for various $E,F\in\LF(\XXX)$.
\eproof

\bRe
If $\XXX=[Y/G]$ as in Example \ref{Ex:[Y/G]-wna} where $k$ is noetherian and $Y=\Spec A$, the finite generation of $H^*(\XXX,\OOO_X)=H^*(G,A)$ over $k$ is proved in \cite{Kallen:Friedlander-Suslin-Thm} by another reduction to \cite[Theorem 10.5]{Kallen:Grosshans-family}. The case $X=BG$ for a finite group scheme $G$ over a field was proved in \cite{Friedlander-Suslin:Cohmology}, and in \cite{Evens:Ring} when $G$ is an abstract finite group.
\eRe

\bRe
If $\XXX$ is of finite type over $\ZZ$, which is sufficient for the proof of Theorem \ref{Th:Main}, instead of \cite[Theorem 10.5]{Kallen:Grosshans-family} one can invoke the earlier \cite[Proposition 55]{Franjou-Kallen:Power-Reductivity}.
\eRe

\section{Residual stacks and gerbes}

In this section, we explicate the residual gerbes of a well-nigh affine stack.

\subsection{Residual stacks}

We recall that algebraic stacks are assumed to be quasi-separated.
For reference we note the following fact.

\bLe
\label{Le:residual}
For an algebraic stack $\XXX$ the following are equivalent.
\begin{enumerate}
\item
\label{It:red-singleton}
$\XXX$ is reduced, and the topological space $|\XXX|$ is a singleton.
\item
\label{It:red-singleton-loc-noeth}
$\XXX$ is reduced and locally noetherian, and $|\XXX|$ is a singleton.
\item
\label{It:SpecL-X}
There is a surjective flat morphism $\Spec L\to\XXX$ where $L$ is a field.
\item
\label{It:SpecL-X-fin-pres}
There is a surjective flat morphism $\Spec L\to\XXX$ locally of finite presentation where $L$ is a field.
\item
\label{It:gerbe/field}
$\XXX$ is a gerbe over a field $\kappa_\XXX$.
\end{enumerate}
The morphism $\XXX\to\Spec \kappa_\XXX$ of \eqref{It:gerbe/field} is unique. For $\Spec L\to\XXX$ as in \eqref{It:SpecL-X-fin-pres} the field extension $\kappa_\XXX\to L$ is finite.
\eLe

\bproof
We have \eqref{It:red-singleton}$\Leftrightarrow$\eqref{It:SpecL-X} by \cite[\href{https://stacks.math.columbia.edu/tag/06MN}{Tag 06MN}]{stacks-project} and \eqref{It:red-singleton-loc-noeth}$\Leftrightarrow$\eqref{It:SpecL-X-fin-pres} by \cite[\href{https://stacks.math.columbia.edu/tag/06MP}{Tag 06MP}]{stacks-project}. Assume that \eqref{It:SpecL-X} holds, and let $\Spec A\to\XXX$ be smooth and surjective. The algebraic space $Z=\Spec A\times_{\XXX}\Spec L$ is quasi-compact and locally of finite type over $L$, thus $Z$ is noetherian. The morphism $Z\to\Spec A$ is flat and surjective, so $A$ is noetherian, thus $\XXX$ is locally noetherian, which implies \eqref{It:red-singleton-loc-noeth}. Assume that \eqref{It:red-singleton-loc-noeth} holds. If $\XXX$ is an algebraic space, then $\XXX$ is the spectrum of a field by \cite[\href{https://stacks.math.columbia.edu/tag/06NH}{Tag 06NH}]{stacks-project}. In general it follows that $\XXX$ is a gerbe over a field by \cite[\href{https://stacks.math.columbia.edu/tag/06QK}{Tag 06QK}]{stacks-project}, which gives \eqref{It:gerbe/field}. Conversely, \eqref{It:gerbe/field} implies \eqref{It:red-singleton-loc-noeth} since a gerbe is smooth over its base by \cite[\href{https://stacks.math.columbia.edu/tag/06QH}{Tag 06QH}, \href{https://stacks.math.columbia.edu/tag/0DLS}{Tag 0DLS}, \href{https://stacks.math.columbia.edu/tag/0DN7}{Tag 0DN7}]{stacks-project}.

The morphism $\XXX\to\Spec\kappa_\XXX$ is unique by \cite[\href{https://stacks.math.columbia.edu/tag/06QD}{Tag 06QD}]{stacks-project}.
If $\Spec L\to\XXX$ is locally of finite presentation, the same holds for $\Spec L\to\Spec\kappa_\XXX$, thus the field extension is finite.
\eproof

If the equivalent conditions of Lemma \ref{Le:residual} hold, $\XXX$ will be called residual. 

\bLe
\label{Le:wna-residual-stack}
For an algebraic stack $\XXX$ the following are equivalent.
\begin{enumerate}
\item
\label{It:residual+wna}
$\XXX$ is residual and well-nigh affine.
\item
\label{It:SpecL-X-finite}
There is a surjective finite locally free morphism $\Spec L\to\XXX$ where $L$ is a field. 
\item
\label{It:gerbe/field-finite}
$\XXX$ is a gerbe with finite inertia over a field $\kappa_\XXX$.
\end{enumerate}
If this holds, $\Spec\kappa_\XXX$ is the coarse moduli space of $\XXX$.
\eLe

\bproof
The implications $\eqref{It:SpecL-X-finite}\Rightarrow\eqref{It:residual+wna}\Rightarrow\eqref{It:gerbe/field-finite}$ are clear. Assume that \eqref{It:gerbe/field-finite} holds, and let $\pi\colon\Spec L\to\XXX$ be surjective flat and locally of finite presentation where $L$ is a field. 
Since a gerbe with finite inertia has finite diagonal
and since $L$ is finite over $\kappa_\XXX$, the morphism $\pi$ is finite, thus finite locally free.
The morphism $\XXX\to\Spec\kappa_\XXX$ has the universal property of the coarse moduli space by \cite[\href{https://stacks.math.columbia.edu/tag/06QD}{Tag 06QD}]{stacks-project}.
\eproof

\subsection{Residual gerbes}

Let $\XXX$ be an algebraic stack (quasi-separated). For $\xi\in|\XXX|$ there is a unique residual substack $\XXX_\xi\to\XXX$ with topological image $\xi$ by \cite[Appendix B]{Rydh:Etale-devissage} or \cite[\href{https://stacks.math.columbia.edu/tag/06RD}{Tag 06RD}]{stacks-project}. The stack $\XXX_\xi$ is a gerbe over a well-defined field $\kappa(\xi)$, and $\XXX_\xi$ is called the residual gerbe of $\XXX$ at $\xi$. 

\bRe
\label{Re:residual-gerbe-functorial}
A morphism of algebraic stacks $f\colon\XXX\to\YYY$ induces a morphism of residual gerbes $f_\xi\colon\XXX_\xi\to\YYY_\eta$ with $\eta=f(\xi)$. If $f$ is a monomorphism, then $f_\xi$ is an isomorphism by the uniqueness of the residual gerbe $\YYY_{\eta}$.
\eRe

\bDe
\label{De:residual-gerbe-wna}
Assume that $\XXX$ is a well-nigh affine stack with coarse moduli space $X=\Spec C$. Then $|\XXX|\cong|X|$, and for $\Fp\in\Spec(C)$ with unique inverse image $\xi_\Fp\in|\XXX|$ we set $\XXX(\Fp)=\XXX_{\xi_\Fp}$, the residual gerbe of $\XXX$ at $\xi_\Fp$.
\eDe

\bRe
The stack $\XXX(\Fp)$ in Definition \ref{De:residual-gerbe-wna} is a gerbe with finite inertia over the field $\kappa(\xi_\Fp)$, in particular $\XXX(\Fp)$ is well-nigh affine by Lemma \ref{Le:wna-residual-stack}. Explicitly, 
\beqn
\label{Eq:XXXFp-first}
\XXX(\Fp)=(\XXX\times_X\Spec k(\Fp))_{\red},
\eeqn
the maximal reduced substack of the fiber of $\XXX$ over $k(\Fp)$, because the latter is reduced and topologically a singleton mapping to $\Fp$ in $X$.  In particular, the morphism $\XXX(\Fp)\to\XXX$ is affine, which again implies that $\XXX(\Fp)$ is well-nigh affine. 
If $\XXX=[Z/{\GL_d}]$ with affine $Z$ as in Lemma \ref{Le:XXX=Z/GLd}, then \eqref{Eq:XXXFp-first} yields
\beqn
\label{Eq:XXXFp-second}
\XXX(\Fp)=[(Z\times_{\XXX}\XXX(\Fp))/{\GL_d}]=[(Z\times_X\Spec k(\Fp))_{\red}/{\GL_d}]
\eeqn
since $Z\times_\XXX\XXX(\Fp)$ is reduced as $Z\to\XXX$ is smooth.
We could take \eqref{Eq:XXXFp-first} or \eqref{Eq:XXXFp-second} as a definition of $\XXX(\Fp)$, but we will need that this stack is a gerbe over a field as a consequence of Lemma \ref{Le:residual}.
\eRe

\bLe
\label{Le:kFp-kappa-insep}
In Definition \ref{De:residual-gerbe-wna}, the field extension $k(\Fp)\to\kappa(\xi_\Fp)$ is purely inseparable. If $\XXX$ is of finite type over a ring, this field extension is finite.
\eLe

\bproof
Since $\XXX\to X$ is a universal homeomorphism, the same holds for the morphism  $\XXX(\Fp)\to\Spec k(\Fp)$, using \eqref{Eq:XXXFp-first}. It also holds for $\XXX(\Fp)\to\Spec\kappa(\xi_\Fp)$, hence for $\Spec\kappa(\xi_\Fp)\to\Spec k(\Fp)$, so this field extension is purely inseparable. Let $Y\to\XXX$ be surjective finite locally with affine $Y$, and choose $y\in Y$ over $\Fp\in X$. We have field extensions $k(x)\to\kappa(\xi_\Fp)\to k(y)$.
If $\XXX$ is of finite type over a ring, then $Y\to X$ is of finite type and integral, hence $k(x)\to k(y)$ is finite, so
$k(x)\to\kappa(\xi_\Fp)$ is finite.
\eproof

\section{Comparison map in the residual case}
\label{Se:Comparison-residual}

In this section we note the following consequence of \cite{Friedlander-Pevtsova:Pi-Supports} or \cite{BIKP:Stratification}; this is Lemma \ref{Le:Intro-rho-XXX-Fp}. The argument is an adaption of \cite[\S7]{Lau:Balmer}.

\bTh
\label{Th:residual-case}
For every residual well-nigh affine stack $\XXX$ the comparison map $\rho_\XXX\colon\Spc(\Perf(\XXX))\to\Spec^h(R_\XXX)$ of \eqref{Eq:rhoX-Intro} is a homeomorphism.
\eTh

\bproof
By Lemma \ref{Le:wna-residual-stack}, $\XXX$ is a gerbe with finite band over a field $K=\kappa_\XXX$. In particular, $\XXX$ is of finite type over $K$, so $\Perf(\XXX)$ is noetherian by Proposition \ref{Pr:R-finite}, and it suffices to show that $\rho_\XXX$ is bijective; see \cite[Corollary 2.8]{Lau:Balmer}.

Assume first that $\XXX$ is a neutral gerbe, thus $\XXX=BG$ for a finite group scheme $G$ over $K$. Let $\modcat(G)$ be the category of finite dimensional representations of $G$ over $K$, thus $\modcat(G)=\Rep_K^{\proj}(G)$ as defined in \S\ref{Se:[Y/G]-gen}. Then $\LF(\XXX)\cong\modcat(G)$ and hence $\Perf(X)\cong D^b(\modcat(G))$ by Proposition \ref{Pr:DbLF-DQcoh-perf}.
In this case, $\rho_\XXX$ is a homeomorphism as a consequence of \cite[Corollary 10.6]{BIKP:Stratification} together with \cite[Remark 2.20]{BCHNP:Quillen}. 

In the general case, $\XXX$ is presented by a finite locally free groupoid of affine schemes $Y_1\rightrightarrows Y_0$ where $Y_0=\Spec L$ for a field $L$, by Lemma \ref{Le:wna-residual-stack}.
The extension $L/K$ is finite by Lemma \ref{Le:residual}. 
Let $K'/K$ be a splitting field of $L/K$ and let $\Gamma=\Aut_K(K')$.
Base change gives a diagram with Cartesian squares
\[
\xymatrix@M+0.2em{
Y_1' \ar[d] \ar@<0.5ex>[r] \ar@<-0.5ex>[r] & Y_0' \ar[d] \ar[r] & \XXX' \ar[d]^f \ar[r] & \Spec K' \ar[d]^g \\
Y_1 \ar@<0.5ex>[r] \ar@<-0.5ex>[r] & Y_0 \ar[r] & \XXX \ar[r]^-q & \Spec K
}
\]
with a compatible action of $\Gamma$ on all schemes and stacks, which is trivial on the lower line. The morphism $f$, or equivalently the left part of the diagram, induces a commutative diagram of topological spaces 
\[
\xymatrix@M+0.2em{
\Spc(\Perf(\XXX')) \ar[r]^-{\rho'} \ar[d]_v & \Spec^h(R_{\XXX'}) \ar[d]^u \\
\Spc(\Perf(\XXX)) \ar[r]^-\rho & \Spec^h(R_\XXX)
}
\]
with a compatible action of $\Gamma$, trivial on the lower line. 

The morphism $\XXX'\to\Spec K'$ is a gerbe by base change, and a neutral gerbe since it has a section induced by $\Spec K'\to\Spec L=Y_0\to\XXX$. Hence $\rho'$ is a homeomorphism by the first case.

The cohomology rings $R=R_\XXX$ and $R'=R_{\XXX'}$ are related by  $R'=R\otimes_KK'$ by \v Cech cohomology (see also Lemma \ref{Le:flat-prop-H}), and $u$ induces a homeomorphism
\[
\bar u\colon\Spec^h(R')/\Gamma \cong \Spec^h(R).
\]
Indeed, let $K''=K'^\Gamma$ and $R''=R\otimes_KK''$. Then $\Gamma$ acts on the ring $R'$ with invariants $R'^\Gamma=R''$, and the field extension $K\to K''$ is purely inseparable, thus a universal homeomorphism on spectra, hence $\Spec^h(R')/\Gamma \cong \Spec ^h(R'') \cong \Spec^h(R)$ by \cite[Lemmas 6.1 and 6.2]{Lau:Balmer}.

Since $\rho'$ and $\bar u$ are homeomorphisms, the same holds for the composition
\[
\Spc(\Perf(\XXX'))/\Gamma \xrightarrow{\;\bar v\;} \Spc(\Perf(\XXX)) \xrightarrow{\;\rho\;} \Spec^h(R).
\]
We will show that $v$ is surjective; this implies that $\rho$ is bijective as required.

The functor $f^*\colon\Perf(\XXX)\to\Perf(\XXX')$ has an exact right adjoint $f_*$ since the morphism $f$ is finite locally free; see Lemma \ref{Le:finite-duality}. By \cite[Theorem 1.7]{Balmer:Surjectivity} it follows that the image of the map $v=\Spc(f^*)$ is equal to $\supp(f_*\OOO_{\XXX'})$. But 
$f_*\OOO_{\XXX'}=q^*g_*\OOO_{\Spec K'}\cong\OOO_{\XXX}^m$ with $m={[K':K]}$, whose support is everything.
\eproof

\section{Fibers of cohomology rings}
\label{Se:Fiber-cohomology}

In this section we prove Lemma \ref{Le:Intro-fiber-Spec-h} along with some extensions.

\subsection{Universal homeomorphisms}

A homomorphism of graded-commu\-tative rings $R\to R'$ will be called a relative universal homeomorphism if the resulting morphism 
\beqn
\label{Eq:rel-univ-homeo}
\Spec R'_{\even}\to\Spec (R_{\even}\otimes_{R_0}R'_0)
\eeqn
of affine schemes is a universal homeomorphism. 

\bRe
\label{Re:univ-homeo}
Let $S\to S'$ be a homomorphism of graded-commutative rings. If the associated map $\Spec(S'_{\even})\to\Spec(S_{\even})$ is a homeo\-morphism, the same holds for the map $\Spec^h(S')\to \Spec^h(S)$ by \cite[Lemma 6.1]{Lau:Balmer}.
\eRe

\subsection{The property H}

For the rest of \S\ref{Se:Fiber-cohomology} let $\XXX$ and $\XXX'$ be well-nigh affine stacks.
We consider the graded-commutative ring
\[
R=R_\XXX=H^*(\XXX,\OOO_\XXX)=\End^*_{\Perf(\XXX)}(\OOO_\XXX).
\]
Here $X=\Spec R_0$ is the coarse moduli space of $\XXX$ as in \S\ref{Se:Coarse} since $R_0=C$.

\bDe
A morphism $\XXX'\to\XXX$ has property H if the resulting homomorphism $R_\XXX\to R_{\XXX'}$ is a relative universal homeomorphism.
\eDe

\bRe
\label{Re:prop-H-compos}
Property H is stable under composition.
\eRe

\bLe
\label{Le:flat-prop-H}
If\/ $\XXX'=\XXX\times_{\Spec k}\Spec k'$ for a flat ring homomorphism $k\to k'$ and a given morphism $\XXX\to \Spec k$, then $R_{\XXX'}=R_\XXX\otimes_kk'$, in particular the natural morphism $\XXX'\to\XXX$ has property H.
\eLe

\bproof
Let $Y_0\to\XXX$ be surjective finite locally free with affine $Y_0$. The associated \v Cech complex $\check C^*=\Gamma(Y_\bullet,\OOO_{Y_\bullet})$ calculates $R=H^*(\check C^*)$. Let $Y_0'=Y_0\times_{\Spec k}\Spec k'$ and $\check C'^*=\Gamma(Y'_\bullet,\OOO_{Y'_\bullet})$. Then $\check C'^*=\check C^*\otimes_kk'$, hence $R'=R\otimes_kk'$ by flatness, and the assertion follows.
\eproof

\bPr
\label{Pr:nilpotent-prop-H}
Let $f\colon\XXX'\to\XXX$ be a closed immersion defined by a nilpotent quasi-coherent ideal $N\!$ of $\OOO_{\XXX}$. Then $f$ has property H.
\ePr

\bproof
Cf.\ \cite[Proposition 8.14]{Lau:Balmer}.
We can assume that $N$ has square zero.
Let $\pi\colon Y\to\XXX$ be surjective finite locally free of degree $d$ with affine $Y$; see Remark \ref{Re:degree}.
By passing to an open covering of $\Spec\ZZ$ we can assume that $d$ is a unit multiple of $p^r$ for a prime $p$. Then $p^r$ annihilates $H^{>0}(\XXX,N)$ by Corollary \ref{Co:HnXN-d}. Let $R=R_\XXX$ and $R'=R_{\XXX'}$.

The exact sequence $0\to N\to\OOO_\XXX\to f_*\OOO_\XXX\to 0$ in $\Qcoh(\XXX)$ induces an exact sequence of graded abelian groups
\beqn
\label{Eq:HXN-R-R'-HXN}
H^*(\XXX,N)\xrightarrow{\;j\;} R\xrightarrow{\;\alpha\;} R'\xrightarrow{\;\delta\;} H^{*+1}(\XXX,N)
\eeqn
where $\alpha$ is a homomorphism of graded-commutative rings, $j$ is a homomorphism of graded $R$-modules whose image is an ideal of square zero, and $\delta$ is a graded derivation. 

Indeed, let $P\to\OOO_\XXX$ be a resolution by projectives in $\Qcoh(\XXX)$ of finite type, using Lemma \ref{Le:D-C}, and let $E=\End^*(P)$ as a dg algebra. Then $H^*(E)=R$ as a graded ring. Let $P'=f^*P$ and $E'=\End^*(P')$. Here $f^*$ preserves projectives since  $f_*$ is exact as $f$ is affine.
Thus $P'\to\OOO_{\XXX'}$ is a resolution by projectives in $\Qcoh(\XXX')$ of finite type, and $R'=H^*(E')$ as a graded ring. We have an exact sequence of complexes \[
0\to P\otimes_{\OOO_\XXX}N\to P\to f_*f^*P\to 0,
\] 
and $\Hom^*(P,-)$ gives an exact sequence $0\to J\to E\xrightarrow\pi E'\to 0$ where $\pi$ is a homomorphism of dg algebras and $J$ is an ideal of square zero. This sequence gives \eqref{Eq:HXN-R-R'-HXN} on the cohomology, and the stated properties of $\alpha$, $j$, and $\delta$ follow.

For a homogeneous element $a\in R'$ of even degree we obtain $\delta(a^{p^r})=0$ and $\delta(p^ra)=0$, thus $a^{p^r}$ and $p^ra$ lie in the image of $\alpha$,
moreover the kernel of $\alpha$ has square zero. Hence the homomorphisms $R_{\even}\to R'_{even}$ and $R_0\to R'_0$ induce universal homeomorphisms on spectra by \cite[\href{https://stacks.math.columbia.edu/tag/0BRA}{Tag 0BRA}]{stacks-project}, and the same follows for $R_{\even}\otimes_{R_0}R_0'\to R_{\even}'$.
\eproof

\subsection{The finite-type noetherian case}

Recall that $\XXX$ is well-nigh affine.

\bPr
\label{Pr:t-prop-H}
Assume that $\XXX$ is of finite type over a noetherian ring, let $t\in \OOO_\XXX(\XXX)$ be an $\OOO_{\XXX}$-regular element, and let $f\colon\XXX'\to\XXX$ be the closed immersion defined by the ideal $t\OOO_\XXX$. Then $f$ has property H.
\ePr

\bproof
Cf.\ \cite[Proposition 8.19]{Lau:Balmer}. Let again $R=R_{\XXX}$ and $R'=R_{\XXX'}$. As in the proof of Proposition \ref{Pr:nilpotent-prop-H} we can assume that $R_{>0}$ is annihilated by $p^r$ for a prime $p$. Proposition \ref{Pr:nilpotent-prop-H} and Remark \ref{Re:prop-H-compos} allow to replace $t$ by $t^m$ for any $m\ge 1$. Hence, since the ring $R$ is noetherian by Proposition \ref{Pr:R-finite}, we can assume that $\Ann_R(t)=\Ann_R(t^2)$.
Let $g\colon\XXX''\to\XXX$ be the closed immersion defined by $t^2$.
There is a diagram of $\OOO_{\XXX}$-modules with exact rows
\[
\xymatrix@M+0.2em{
0 \ar[r] & \OOO_\XXX \ar[r]^t \ar[d] & \OOO_\XXX \ar[r] \ar[d] & f_*\OOO_{\XXX'} \ar[r] \ar[d]^{\id} & 0 \\
0 \ar[r] & f_*\OOO_{\XXX'} \ar[r]^{\bar t} & g_*\OOO_{\XXX''} \ar[r] & f_*\OOO_{\XXX'} \ar[r] & 0,\!
}
\]
which induces a commutative diagram of graded abelian groups with exact upper row, where $\delta$ and $\delta'$ have degree $1$,
\[
\xymatrix@M+0.2em{
R \ar[r]^t & R \ar[r]^\pi & R' \ar[r]^\delta \ar[d]_{\id} & R \ar[r]^t \ar[d]^{\pi} & R \\
& & R'\ar[r]^{\delta'} & R'.
}
\]
For an an element $a\in R'$ of even degree we have $\pi\delta(a^{p^r})=0$ since $\delta'$ is a graded derivation by the proof of Proposition \ref{Pr:nilpotent-prop-H}, hence $\delta(a^{p^r})=tb$ for some $b\in R$, then $t^2b=t\delta(a^{p^r})=0$, and $tb=0$ by the annihilator equality. Moreover $\delta(p^ra)=0$. Hence $a^{p^r}$ and $p^ra$ lie in the subring $R/tR\subseteq R'$. As in the proof of Proposition \ref{Pr:nilpotent-prop-H} it follows that $R/tR\to R'$ is a relative universal homeomorphism. The same holds for $R\to R/tR$ and hence for $R\to R'$.
\eproof

As earlier let $\XXX\to X=\Spec C$ be the coarse moduli space (\S\ref{Se:Coarse}).

\bLe
\label{Le:find-property-H}
Assume that $\XXX$ is of finite type over a noetherian ring $k$. Let $\Fp\in X$ be a maximal ideal. If the morphism $\XXX(\Fp)\to\XXX$ of Definition \ref{De:residual-gerbe-wna} is not an isomorphism, then it factors into closed immersions $\XXX(\Fp)\to\XXX'\xrightarrow u\XXX$ where $u$ has property H and $u$ is not an isomorphism.
\eLe

\bproof
Cf.\ \cite[Proposition 8.23]{Lau:Balmer}.
The morphism $\XXX(\Fp)\to\XXX$ is a closed immersion by \eqref{Eq:XXXFp-first}. If $\XXX$ is not reduced, we can take $\XXX'=\XXX_{\red}$, using Proposition \ref{Pr:nilpotent-prop-H}. Assume that $\XXX$ is reduced. Since $\XXX(\Fp)\to\XXX$ is not an isomorphism, $\XXX$ is not residual, so $|\XXX|=|X|$ has more than one point; see Lemma \ref{Le:residual}. 
Let $\XXX=[Z/{\GL_d}]$ with $Z=\Spec S$ as in Lemma \ref{Le:XXX=Z/GLd}. Then $S$ is reduced and of finite type over $k$, in particular the ring $S$ is noetherian. We recall that $X=\Spec C$ with $C=H^0(\XXX,\OOO_\XXX)=H^0(\GL_d,S)$.

Assume that a minimal prime $\Fq$ of $S$ maps to $\Fp$ in $X$. Let $V\subseteq Z$ be an open neighborhood of $\Fp$ that avoids all other minimal primes of $S$. The image of $V\to\XXX$ is an open substack $\UUU\subseteq\XXX$ which lies over $\Fp\in X$. Hence $\UUU$ is residual and thus $\UUU=\XXX(\Fp)$ by the uniqueness of residual gerbes. Thus $X(\Fp)\subseteq\XXX$ is open and closed, hence $\{\Fp\}\subseteq X$ is open and closed. Since $C$ is reduced this yields $C=\Spec k(\Fp)\times C_1$ as rings, so $C\to k(\Fp)$ is flat, and we can take $\XXX'=\XXX\times_X\Spec k(\Fp)$, using Lemma \ref{Le:flat-prop-H}.

Assume that no minimal prime of $S$ maps to $\Fp$ in $X$. By prime avoidance there is an element $t\in\Fp$ that avoids all minimal primes of $S$. Since $S$ is reduced, $t$ is an $S$-regular element, thus $t$ is $\OOO_\XXX$-regular, und we can take the closed substack $\XXX'\subseteq\XXX$ defined by $t$, using Proposition \ref{Pr:t-prop-H}.
\eproof

\bPr
\label{Pr:residual-stack-property-H-noeth}
Assume that $\XXX$ is of finite type over a noetherian ring. For any $\Fp\in X$, the morphism $\XXX(\Fp)\to\XXX$ has property H.
\ePr

\bproof
Cf.\ \cite[Proposition 8.24]{Lau:Balmer}.
We can assume that $\Fp$ is maximal because our morphism factors as
$\XXX(\Fp)\to\XXX\times_X\Spec C_\Fp\to\XXX$ where the second morphism has property H by Lemma \ref{Le:flat-prop-H}. Since $\XXX$ is noetherian, there is a minimal closed substack $\XXX_1\subseteq\XXX$ with $\XXX(\Fp)\subseteq\XXX_1$ such that $\XXX_1\to\XXX$ has property H. We have $\XXX(\Fp)=\XXX_1(\Fp)$ by Remark \ref{Re:residual-gerbe-functorial} or by \eqref{Eq:XXXFp-first}.
By Lemma \ref{Le:find-property-H}, the minimality of $\XXX_1$ implies that $\XXX(\Fp)=\XXX_1$.
\eproof

\bCo
\label{Co:residual-stack-property-H}
Assume that $\XXX$ is of finite type over a noetherian ring. 
For $\Fp\in X$, the homomorphism $R_\XXX\to R_{\XXX(\Fp)}$ induces a homeomorphism between $\Spec^h(R_{\XXX(\Fp)})$ and the fiber over $\Fp$ of the map $\Spec^h(R_\XXX)\to\Spec(C)$.
\eCo

\bproof
Since the field extension $k(\Fp)\to\kappa(\xi_\Fp)$ is purely inseparable by Lemma \ref{Le:kFp-kappa-insep}, the said fiber is homeomorphic to $\Spec^h(R_\XXX\otimes_C\kappa(\xi_\Fp))$, which is homeomorphic to $\Spec^h(R_{\XXX(\Fp)})$ by Proposition \ref{Pr:residual-stack-property-H-noeth} (see also Remark \ref{Re:univ-homeo}).
\eproof

\subsection{Noetherian approximation}

The following additions are not necessary for the proof of the main result. We begin with an extension of Lemma \ref{Le:flat-prop-H}.
Recall that $\XXX$ is well-nigh affine.

\bPr
\label{Pr:bc-prop-H}
If $\XXX'=\XXX\times_{\Spec k}\Spec k'$ for any ring homomorphism $k\to k'$ and a given morphism $\XXX\to \Spec k$, then $\XXX'\to\XXX$ has property H.
\ePr

\bLe
\label{Le:bc-prop-H-finite-type}
In the proof of Proposition \ref{Pr:bc-prop-H} we can assume that $k$ is noetherian and that $\XXX$ and $k'$ are of finite type over $k$.
\eLe

\bproof
We can assume that $k=\Gamma(\XXX,\OOO_\XXX)$.
If a well-nigh affine stack $\YYY$ is a filtered limit of well-nigh affine stacks $\YYY_i$ with affine transition maps (Definition \ref{De:limit-stack-affine}), then $R_{\YYY}=\colim_iR_{\YYY_i}$ by \v Cech cohomology. Hence, if a morphism $f\colon\YYY'\to\YYY$ of well-nigh affine stacks is a filtered limit with affine transition maps of morphisms $f_i\colon\YYY'_i\to\YYY_i$ with property H, then $f$ has property H. Thus  we can assume that $k'$ is a $k$-algebra of finite presentation. Moreover, if $\XXX=\lim_i\XXX_i$ as in Lemma \ref{Le:approximation} and $k_i=\Gamma(\XXX_i,\OOO_{\XXX_i})$, then $k=\colim_ik_i$, the homomorphism $k\to k'$ arises from a finitely presented homomorphism $k_i\to k'_i$ for some $i$, and we can assume that $\XXX=\XXX_i$.
\eproof

\bproof[Proof of Proposition \ref{Pr:bc-prop-H}]
We assume that $k$ is noetherian and that $\XXX$ and $k'$ are of finite type over $k$.
Since $k'$ is a quotient of a polynomial ring over $k$, by Lemma \ref{Le:flat-prop-H} we can assume that $k\to k'$ is surjective. Then $f\colon\XXX'\to\XXX$ is a closed immersion, in particular $f_*\OOO_{\XXX'}$ is coherent, hence $R_{\XXX'}$ is a finite module over $R_{\XXX}$ by Proposition \ref{Pr:R-finite}. It follows that the morphism \eqref{Eq:rel-univ-homeo} with $R=R_\XXX$ and $R'=R_{\XXX'}$ is finite, hence universally closed. Moreover, this morphism is universally bijective since for every $\Fp'\in\Spec(R_0')$ with image $\Fp\in\Spec(R_0)$, the resulting morphism between the fibers over $\Fp'$ is a universal homeomorphism by Corollary \ref{Pr:residual-stack-property-H-noeth}; note that $k(\Fp')$ is a purely inseparable extension of $k(\Fp)$ as a consequence of Lemma \ref{Le:kFp-kappa-insep}.
\eproof

Now we can drop the noetherian hypothesis in Proposition \ref{Pr:residual-stack-property-H-noeth}.

\bPr
\label{Pr:residual-stack-property-H}
For $\Fp\in X$, the morphism $\XXX(\Fp)\to\XXX$ has property H.
\ePr

\bproof
By Proposition \ref{Pr:bc-prop-H} we can replace $\XXX$ by $\XXX\times_{X}\Spec k(\Fp)$. Then $f\colon\XXX(\Fp)\to\XXX$ is a closed immersion defined by a locally nilpotent ideal $N\subseteq\OOO_{\XXX}$ by \eqref{Eq:XXXFp-first}. Let $\XXX=\lim_i\XXX_i$ as in Lemma \ref{Le:approximation}, and let $f_i\colon\ZZZ_i\to\XXX_i$ be the closed substack defined by the inverse image of $N$ in $\OOO_{\XXX_i}$. Then $f_i$ is defined by a nilpotent ideal, hence $f_i$ has property H by Proposition \ref{Pr:nilpotent-prop-H}. Moreover $f=\lim_if_i$. By the proof of Lemma \ref{Le:bc-prop-H-finite-type} it follows that $f$ has property H.
\eproof

\bRe
By replacing Proposition \ref{Pr:residual-stack-property-H-noeth} with Proposition \ref{Pr:residual-stack-property-H}, Corollary \ref{Co:residual-stack-property-H} holds without the noetherian hypothesis. This proves Lemma \ref{Le:Intro-fiber-Spec-h}.
\eRe

\section{Fibers of Balmer spectra}
\label{Se:Fiber-Balmer}

In this section we prove Lemma \ref{Le:Intro-fiber-Spc}.

\subsection{Tensor nilpotence}

A functor $\varphi\colon\TTT\to\TTT'$ of tensor triangulated categories detects tensor nilpotence of morphisms if for $h\colon A\to B$ in $\TTT$ with $\varphi(h)=0$ we have $h^{\otimes n}=0$ for some $n$. If $\varphi$ detects tensor nilpotence of morphisms and $\TTT$ is rigid, then $\Spc(\varphi)\colon\Spc(\TTT')\to\Spc(\TTT)$ is surjective by \cite[Theorem 1.1]{Balmer:Surjectivity}.

\bPr
\label{Pr:detect-tensor-nilpotence}
Let $\XXX$ be a well-nigh affine stack and $f\colon\XXX'\to\XXX$ a morphism such that one of the following holds.
\begin{enumerate}
\item
\label{It:detect-nilpotence-nil}
The morphism $f$ is a closed immersion defined by a nilpotent quasi-coherent ideal $N$
of $\OOO_{\XXX}$.
\item
\label{It:detect-nilpotence-t}
We have $\XXX'=\UUU\sqcup\ZZZ_1$ where $\ZZZ_1$ is a closed substack of $\XXX$ defined by an $\OOO_{\XXX}$-regular element $t\in\Gamma(\XXX,\OOO_\XXX)$, and $\UUU=\XXX\setminus\ZZZ_1$.
\end{enumerate}
Then $f^*\colon\Perf(\XXX)\to\Perf(\XXX')$ detects tensor nilpotence of morphisms.
\ePr

\bproof
Cf.\ \cite[Proposition 9.1]{Lau:Balmer}.
Let $\alpha\colon P\to Q$ in $\Perf(\XXX)$ with $f^*(\alpha)=0$. Assume that $P$ and $Q$ are as in Remark \ref{Re:PQ-special}.
Then $\alpha$ is represented by a homomorphism of complexes $u\colon P\to Q$ which is unique up to homotopy.

In the case \eqref{It:detect-nilpotence-nil} we have an exact sequence $0\to Q\otimes N\to Q\xrightarrow\varepsilon f_*f^*Q\to 0$ of complexes in $\Qcoh(\XXX)$, and $f^*(\alpha)=0$ means that $P\xrightarrow uQ\xrightarrow\varepsilon f_*f^*Q$ is homotopic to zero. Since the corresponding homotopy can be lifted over $\varepsilon$, we can assume that $u\colon P\to Q\otimes N$. If $N^r=0$ if follows that $u^{\otimes r}=0$.

In the case \eqref{It:detect-nilpotence-t} let $\ZZZ_r$ be the closed substack of $\XXX$ defined by $t^r$, and let $g\colon\UUU\to\XXX$ and $h_r\colon\ZZZ_r\to\XXX$ be the natural morphisms. Let $C=\Gamma(\XXX,\OOO_\XXX)$. Remark \ref{Re:PQ-special} implies that 
\[
\Hom_{\!\:\Perf(\UUU)}(g^*P,g^*Q)=\Hom_{\!\:\Perf(\XXX)}(P,Q)\otimes_CC_t.
\]
Hence $g^*(\alpha)=0$ gives $t^r\alpha=0$ for some $r$, and the exact triangle
\[
P\xrightarrow {t^r}P\xrightarrow\varepsilon h_{r*}h_r^*P\to P[1]
\]
in $D(\Qcoh(\XXX))$ implies that $\alpha$ factors as $P\xrightarrow\varepsilon h_{r*}h_r^*P\xrightarrow{\bar\alpha}Q$ in $D(\Qcoh(\XXX))$. Now we have the following commutative diagram in $D(\Qcoh(\XXX))$.
\[
\xymatrix@M+0.2em{
P^{\otimes r}\otimes P \ar[r]^-{1\otimes\varepsilon} \ar[d]^{\alpha^{\otimes r}\otimes 1} &
P^{\otimes r}\otimes h_{r*}h_r^*P \ar[r]^-\sim \ar[d]^{\alpha^{\otimes r}\otimes 1} &
h_{r*}(h_r^*P^{\otimes r}\otimes h_r^*P) \ar[d]^{h_{r*}(h_r^*\alpha^{\otimes r}\otimes 1)} \\
Q^{\otimes r}\otimes P \ar[r]^-{1\otimes\varepsilon} &
Q^{\otimes r}\otimes h_{r*}h_r^*P \ar[r]^-\sim &
h_{r*}(h_r^*Q^{\otimes r}\otimes h_r^*P)
}
\]
We have $h_1^*\alpha=0$ and hence $h_r^*\alpha^{\otimes r}=0$ by the proof of \eqref{It:detect-nilpotence-nil}. It follows that $\alpha^{\otimes r+1}=(1\otimes\bar\alpha)\circ(1\otimes\varepsilon)\circ(\alpha^{\otimes r}\otimes 1)$ is zero.
\eproof

\subsection{The property S}

Let $\XXX$ be a well-nigh affine stack with coarse moduli space $X=\Spec C$ as in \S\ref{Se:Coarse}, thus $C=\Gamma(\XXX,\OOO_\XXX)=\End(\OOO_\XXX)$. 
There is a natural continuous map 
\cite[Corollary 5.6]{Balmer:Spectra3}
\beqn
\label{Eq:bar-rho-XXX}
\bar\rho_\XXX\colon\Spc(\Perf(\XXX))\to\Spec(C)\cong|\XXX|.
\eeqn
For a morphism of well-nigh affine stacks $f\colon\XXX'\to\XXX$, by functoriality with respect to $f^*\colon\Perf(\XXX)\to\Perf(\XXX')$ we obtain a continuous map
\beqn
\label{Eq:f*-Spc-rel}
f_*^{\rel}\colon \Spc(\Perf(\XXX'))\to\Spc(\Perf(\XXX))\times_{|\XXX|}|\XXX'|.
\eeqn

\bDe
$f\colon\XXX'\to\XXX$  has property S if $f_*^{\rel}$ is surjective.
\eDe

\bRe
The property S is stable under composition.
\eRe

\bLe
\label{Le:flat-prop-S}
If $\XXX'=\XXX\times_{\Spec k}\Spec k'$ for a localisation homomorphism $k\to k'$ and a given morphism $\XXX\to \Spec k$, then $\XXX'\to\XXX$ has property S.
\eLe

\bproof
Cf.\ \cite[Proposition 10.6]{Lau:Balmer}.
We can assume that $X=\Spec k$ is the coarse moduli space of $\XXX$,
and then $X'=\Spec k'$ is the coarse moduli space of $\XXX'$ since $k\to k'$ is flat. Let $k'=S^{-1}k$. Remark \ref{Re:PQ-special} implies that the obvious functor $\varphi\colon S^{-1}\Perf(\XXX)\to\Perf(\XXX')$ is fully faithful. We obtain a commutative diagram of topological spaces
\[
\xymatrix@M+0.2em{
\Spc(\Perf(\XXX')) \ar[r]^-{\Spc(\varphi)} \ar[d] &
\Spc(S^{-1}\Perf(\XXX)) \ar[r] \ar[d] &
\Spc(\Perf(\XXX)) \ar[d] \\
\Spec(k') \ar[r]^{\id} &
\Spec(k') \ar[r] &
\Spec(k)
}
\]
where the second square is a fiber product by \cite[Theorem 5.4]{Balmer:Spectra3}, and the map $\Spc(\varphi)$ is surjective since $\varphi$ is fully faithful, using \cite[Corollary 1.8]{Balmer:Surjectivity}.
\eproof

\bPr
\label{Pr:nilpotent-prop-S}
If $f\colon\XXX'\to\XXX$ is a closed immersion defined by a nilpotent ideal, then $f$ has property S.
\ePr

\bproof
Cf.\ \cite[Proposition 10.8]{Lau:Balmer}.
Since $\XXX'\to\XXX$ is a homeomorphism we have to show that $\Spc(\Perf(\XXX'))\to\Spc(\Perf(\XXX))$ is surjective, which follows from Proposition \ref{Pr:detect-tensor-nilpotence} and \cite[Theorem 1.1]{Balmer:Surjectivity}.
\eproof

\bPr
\label{Pr:t-prop-S}
If $f\colon\XXX'\to\XXX$ is a closed immersion defined by an $\OOO_{\XXX}$-regular element $t\in\Gamma(\XXX,\OOO_\XXX)$, then $f$ has property S.
\ePr

\bproof
Cf.\ \cite[Proposition 10.10]{Lau:Balmer}.
Let $C=\Gamma(\XXX,\OOO_\XXX)$.
We consider the stacks $\XXX_i=\XXX\times_{\Spec C}\Spec C_i$ for $C_1=C_t$, $C_2=C/tC$, and $C_3=C_1\times C_2$; thus $\XXX'=\XXX_2$ and $\XXX_3=\XXX_1\sqcup\XXX_2$. The map $|\XXX_3|\to|\XXX|$ is bijective since $\XXX\to\Spec C$ is a universal homeomorphism and $\Spec(C_3)\to\Spec(C)$ is bijective. The map $\Spc(\Perf(\XXX_3))\to\Spc(\Perf(\XXX))$ is surjective by Proposition \ref{Pr:detect-tensor-nilpotence} and \cite[Theorem 1.1]{Balmer:Surjectivity}. Thus the morphism $\XXX_3\to\XXX$ has property S. This carries over to $\XXX_2\to\XXX$ because the decomposition $\XXX_3=\XXX_1\sqcup\XXX_2$ induces a similar decomposition of $\Spc(\Perf(-))$.
\eproof

\bLe
\label{Le:find-property-S}
Lemma \ref{Le:find-property-H} holds with property S in place of property H.
\eLe

\bproof
In the proof of Lemma \ref{Le:find-property-H}, replace 
Lemma \ref{Le:flat-prop-H} and
Propositions \ref{Pr:nilpotent-prop-H} and \ref{Pr:t-prop-H}
by 
Lemma \ref{Le:flat-prop-S} and
Propositions \ref{Pr:nilpotent-prop-S} and \ref{Pr:t-prop-S}.
\eproof

\bPr
\label{Pr:residual-stack-property-S}
Assume that $\XXX$ is of finite type over a noetherian ring. For any $\Fp\in X$, the morphism $\XXX(\Fp)\to\XXX$ has property S.
\ePr

\bproof
In the proof of Proposition \ref{Pr:residual-stack-property-H-noeth}
replace
Lemma \ref {Le:flat-prop-H} and
Lemma \ref{Le:find-property-H}
by
Lemma \ref {Le:flat-prop-S} and
Lemma \ref{Le:find-property-S}.
\eproof

\bRe
\label{Re:residual-stack-property-S}
Proposition \ref{Pr:residual-stack-property-S} means that the natural map $\Spc(\Perf(\XXX(\Fp)))\to\Spc(\Perf(\XXX))$ has image equal to the fiber $(\bar\rho_\XXX)^{-1}(\Fp)$. This is Lemma \ref{Le:Intro-fiber-Spc}.
\eRe

\section{Main result}

\bTh
For every well-nigh affine stack $\XXX$, the comparison map 
\[
\rho_\XXX\colon\Spc(\Perf(\XXX))\to\Spec^h(R_\XXX)
\] 
of \eqref{Eq:rhoX-Intro} is a homeomorphism.
\eTh

\bproof
Let $R=R_\XXX$.
We can assume that $\XXX$ is of finite type over $\ZZ$ because for a filtered limit $\XXX=\lim_i\XXX_i$ as in Lemma \ref{Le:approximation} and $R_i=R_{\XXX_i}$ we have $\Perf(\XXX)\cong\colim_i\Perf(\XXX_i)$ by Proposition \ref{Pr:PerfXXX-colim}, thus $R\cong\colim_i R_i$, which gives homeomorphisms $\Spec^h(R)\cong\lim_i\Spec^h(R_i)$ and $\Spc(\Perf(\XXX))\cong\lim_i\Spc(\Perf(\XXX_i))$. 

Let $X=\Spec C$ be the coarse moduli space of $\XXX$, so $C=\Gamma(\XXX,\OOO_\XXX)=R_0$. By functoriality with respect to the morphism $\XXX\to X$ we have the following commutative diagram of topological spaces.
\beqn
\label{Eq:fund-diag}
\xymatrix@M+0.2em{
\Spc(\Perf(\XXX)) \ar[r]^-{\rho_\XXX} \ar[d] &
\Spec^h(R) \ar[d] \\
\Spc(\Perf(X)) \ar[r]^-{\rho_X} &
\Spec(C)
}
\eeqn

Here $\rho_X$ is a homeomorphism by \cite[Proposition 8.1]{Balmer:Spectra3}, based on \cite[Theorem 3.15]{Thomason:Classification}. The diagonal arrow is the map $\bar\rho_{\XXX}$ of \eqref{Eq:bar-rho-XXX}. The category $\Perf(\XXX)$ is noetherian by Proposition \ref{Pr:R-finite}. Hence the map $\rho_\XXX$ is a homeomorphism iff it is bijective by \cite[Corollary 2.8]{Lau:Balmer}. 
This holds iff for each $\Fp\in\Spec(C)$, the map $\rho_\XXX$ induces a bijective map between the fibers over $\Fp$.

Let $\XXX(\Fp)$ be the residual gerbe of $\XXX$ over $\Fp$ as in Definition \ref{De:residual-gerbe-wna}. By functoriality we obtain a commutative diagram of topological spaces, where an index $\Fp$ means fiber over $\Fp$ in \eqref{Eq:fund-diag}.
\[
\xymatrix@M+0.2em{
\Spc(\Perf(\XXX(\Fp))) \ar[r]^-{\rho_{\XXX(\Fp)}} \ar[d]_u &
\Spec^h(R_{\XXX(\Fp)}) \ar[d]^j \\
\Spc(\Perf(\XXX))_\Fp \ar[r]^-{(\rho_{\XXX})_\Fp} &
\Spec^h(R_\XXX)_\Fp  \\
}
\]

Here $\rho_{\XXX(\Fp)}$ is a homeomorphism by Theorem \ref{Th:residual-case}, $j$ is a homeomorphism by Corollary \ref{Co:residual-stack-property-H}, and $u$ is surjective by Proposition \ref{Pr:residual-stack-property-S} (see Remark \ref{Re:residual-stack-property-S}). Hence $(\rho_\XXX)_\Fp$ is bijective as required.
\eproof


\end{document}